\documentclass[12pt,a4paper]{amsart}
\usepackage[utf8]{inputenc}

\usepackage{amsfonts}
\usepackage{amsthm}
\usepackage{amsmath}
\usepackage{amscd}
\usepackage{t1enc}
\usepackage[mathscr]{eucal}
\usepackage{indentfirst}
\usepackage{graphicx}
\usepackage{graphics}
\usepackage{pict2e}
\usepackage{epic}
\usepackage{color}
\numberwithin{equation}{section}
\usepackage[margin=2.9cm]{geometry}
\usepackage{epstopdf} 
\usepackage{hyperref}

\theoremstyle{plain}
\newtheorem{Th}{Theorem}[section]
\newtheorem{lemma}[Th]{Lemma}
\newtheorem{Cor}[Th]{Corollary}

 \theoremstyle{definition}
\newtheorem{Def}[Th]{Definition}

\newtheorem{?}[Th]{Problem}

\newcommand{\NN}{\mathbb{N}}
\newcommand{\ZZ}{\mathbb{Z}}
\newcommand{\QQ}{\mathbb{Q}}
\newcommand{\RR}{\mathbb{R}}

\newcommand{\lnko}[1]{\text{gcd}(#1)}
\newcommand{\kar}[1]{\chi_t^{(#1)}(m)}
\newcommand{\kkar}[1]{{\widehat{\chi}_t^{(#1)}(m)}}

\newcommand{\scal}[2]{\langle #1,#2\rangle}
\newcommand{\rg}{\widetilde{\gamma}}
\newcommand{\ig}{\widehat{\gamma}}

\begin{document}
\title{Estimating the $p$-adic valuation of the resultant}
\author[K. Sz.]{Krist\'{o}f Szab\'{o}}
\address{ELTE: E\"{o}tv\"{o}s Lor\'{a}nd University \\ H-1117 Budapest
\\ P\'{a}zm\'{a}ny P\'{e}ter s\'{e}t\'{a}ny 1/C}

\thanks{I am particularly grateful for the assistance given by Dr Gergely Zábrádi. In addition, the 2021/2022 UNKP supported this research. }

\begin{abstract}
    Let $f$ and $g$ be two monic polynomials with integer coefficients and nonzero resultant $r$. 
    Assume that $v_p(f(n))\ge s_1$ and $v_p(g(n))\ge s_2$ hold for all integers $n$ for some $s_1, s_2$ fixed non-negative integers. Let $S$ denote the maximum of $v_p(gcd(f(n),g(n)))$ over all integers $n$. In this paper, we establish multiple lower bound for $v_p(r)$. More specifically, we show that $v_p(r)\ge S-\max{s_1,s_2}+ps_1s_2\frac{p-1}{p-p^{-k}}$, where $k=\lfloor \log_p((p-1)\max\{s_1,s_2\}+1)\rfloor -1$.
\end{abstract}
\maketitle

\section{Introduction}

Let $f,g\in\ZZ[x]$ be monic polynomials with resultant $r$. Some properties of the greatest common divisor of $f(n)$ and $g(n)$ was examined by P. Frenkel and J. Pelik\'an in their paper \cite{frenkel2017greatest}. They managed to prove the following theorems:
\begin{enumerate}
    \item For every integer $n\in \ZZ$ the $\lnko{f(n),g(n)}$ divides $r$.
    \item If $r$ square-free, then every divisor of $r$ will occur in the series $(\lnko{f(n),g(n)})_{n\in\ZZ}$.
    \item If $p^p$ does not divide $r$ for every $p$ prime, then there exists an $n\in\ZZ$ integer such as $\lnko{f(n),g(n)}=1$.
\end{enumerate}
As the aforementioned claims show, the series $(\lnko{f(n),g(n)})_{n\in\ZZ}$ and the $r$ resultant are strongly related. This connection motivated further investigation of the topic. Let $v_p$ be the $p$-adic valuation and assume that $v_p(\lnko{f(n),g(n)})=\min{v_p(f(n),v_p(g(n))}\ge s$ for every $n\in\ZZ$. Furthermore, $S$ denotes the maximum, $S=\max_{n\in\ZZ}{\min{v_p(f(n)),v_p(g(n))}}$. This notion helps us to rewrite the first result of P. Frenkel and J. Pelikán to the form $v_p(r)\ge S$. In the paper \cite{frenkel2018estimating} the authors, P. Frenkel and G. Z\'abr\'adi, were focusing on this inequality, and their interest was finding the best possible lower bound for $v_p(r)$ in terms of $S,s$ parameters. A short summary of their paper is presented below:
\begin{enumerate}
    \item The $v_p(r)\ge S-s+B(2s)-2B(s)$ inequality holds for every non-negative $s$, where $B$ is a function defined in their paper.
    \item The corollary of (1) is the following inequality:
    $$v_p(r)\ge S-s+(p-1)s^2$$ This lower bound is asymptotically sharp.
    \item For every $s\le p$ we have $v_p(r)\ge S-s+ps^2$ and it is sharp.
    \item In the case of equality in the last statement, the series $(\min{v_p(f(n)),v_p(g(n))})_{n\in\ZZ}$ takes every value from the interval $[s,S]$.
\end{enumerate}

In this paper we establish a new lower bound for $v_p(r)$ and also examine more general problems. The main results of this paper are the following theorems.
\begin{Th}
Let $f,g\in\ZZ[x]$ be monic polynomials with nonzero $r$ resultant and assume that $\min\{v_p(f(n),v_p(g(n))\}\ge s$ for every $n\in \ZZ$. Then we have the following lower bound for $v_p(r)$.
\begin{align*}
   v_p(r)\ge S-s+ps^2\frac{p-1}{p-p^{1-\lfloor \log_p((p-1)s+1)\rfloor}}
\end{align*}
\end{Th}

\bigskip
\begin{Th}
Let $f,g\in\ZZ[x]$ be monic polynomials with nonzero $r$ resultant and assume that $v_p(f(n))\ge s_1$ and $v_p(g(n))\ge s_1$ for every $n\in \ZZ$. Then we have the following lower bound for $v_p(r)$.
\begin{align*}
    v_p(r)\ge S-\max{s_1,s_2}+ps_1s_2\frac{p-1}{p-p^{-k}}
\end{align*}
Where $k=\lfloor \log_p((p-1)\max\{s_1,s_2\}+1)\rfloor -1$.
\end{Th}

\section{An approach with $p$-adic numbers}
Throughout this chapter $f,g$ will denote monic polynomials with integer coefficients and with nonzero $r$ resultant. Furthermore, assume that $v_p(f(n))\ge s_1$ and $v_p(g(n))\ge s_2$ for every $n\in \ZZ$. (Note: $s_1,s_n\in \NN$) Let $K$ be the splitting field of the polynomials $f,g$ and write $f,g$ as $f(x)=\prod_{i=1}^k(x-\alpha_i)$, $g(x)=\prod_{j=1}^l(x-\beta_j)$. From the definition of $v_p$ it follows that
\begin{align}
v_p(f(n))&=\sum_{i=1}^k v_p(n-\alpha_i)\ge s_1\quad \forall ~n\in\ZZ  \\
v_p(g(n))&=\sum_{j=1}^l v_p(n-\beta_j)\ge s_2 \quad \forall ~n\in\ZZ  \\
v_p(r)&=\sum_{i,j} v_p(\alpha_i-\beta_j)
\end{align}
\begin{Def}
\label{karDef}
For arbitrary $(x-\delta)\in\overline{\QQ_p}[x]$ linear polynomial, $m\in\ZZ_p$ $p$-adic integer and $t\in\RR_{\ge0}$ non-negative real we define
$$\kar{x-\delta}=\begin{cases}1&~~\text{if}~ v_p(m-\delta)\ge t\\0&~~\text{otherwise}\end{cases}$$
In general, for higher degree polynomials we define $\kar{h}$ to be $\kar{h}=\sum_i\kar{x-\delta_i}$, where $h(x)=\prod_i(x-\delta_i)$.
\end{Def}

\begin{lemma}
For arbitrary $h\in\ZZ[x]$ monic polynomial and $m\in\ZZ_p$ $p$-adic integer
\begin{enumerate}
    \item $\kar{h}$ is a monotonically decreasing function in $t$.
    \item $v_p(h(m))=\int_0^\infty \kar{h}dt$
\end{enumerate}
\end{lemma}
\begin{proof}
If $h(x)=x-\delta$ is a linear monic polynomial, then
\begin{align*}
    \kar{x-\delta}=\begin{cases}
    1~~&\text{if}~ t\in [0,v_p(m-\delta)]
    \\0~~&\text{otherwise}
    \end{cases}
\end{align*}
Therefore $\kar{x-\delta}$ is monotonically decreasing and $\int_0^\infty \kar{h}dt=v_p(m-\delta)$.
\\The statements for higher degree polynomials follow from the linear aspect of $\kar{h}$. Let $h\in\overline{\QQ_p}[x]$ be a monic polynomials and $h(x)=\prod_i (x-\delta_i)$ be its factorization. The sum of monotonically decreasing functions is monotonically decreasing and $\kar{h}=\sum_i\kar{x-\delta_i}$, therefore $\kar{h}$ is indeed monotonically decreasing. On the other hand
$$\int_0^\infty \kar{h}dt=\sum_i\int_0^\infty \kar{x-\delta_i}dt=\sum_i v_p(m-\delta_i)=v_p(h(m))$$
which proves (2).
\end{proof}
\begin{lemma}\label{vprchiint}
$$v_p(r)\ge \int_0^\infty \sum_{m=1}^{p^{\lceil t\rceil}}\kar{f}\kar{g} dt$$
\end{lemma}
\begin{proof}
Notice that $v_p(r)=\sum_{i,j}(\alpha_i-\beta_j)$ and $v_p(\alpha_i-\beta_j)\ge \min\{v_p(m-\alpha_i),v_p(m-\beta_j)\}$ for every $m\in\ZZ_p$ due to the triangular inequality. It follows that
$$v_p(r)=\sum_{i,j}v_p(\alpha_i-\beta_j)\ge\sum_{i,j}\max_{m\in\ZZ_p}\min\{v_p(m-\alpha_i),v_p(m-\beta_j)\}$$
We need one last observation to transform the last inequality into a much convenient one. Notice that
$$\min\{v_p(m-\alpha_i),v_p(m-\beta_j)\}=\int_0^\infty \kar{x-\alpha_i}\kar{x-\beta_j}dt$$
The $\kar{x-\alpha_i}\kar{x-\beta_j}$ takes $1$ if and only if both terms are $1$ which is equivalent to saying that $t\in[0,v_p(m-\alpha_i)]\cap[0,v_p(m-\beta_j)]=[0,\min\{v_p(m-\alpha_i),v_p(m-\beta_j)\}]$. Therefore the right hand side is equal to the length of the interval within the integrand takes $1$.
\\The next step is finding a way to effortlessly determine $m\in\ZZ_p$ which maximizes the aformentioned quantity. For given $t$ the $l(m):=\kar{x-\alpha_i}\kar{x-\beta_j}$ takes $1$ for at most one $m$ from the set $\{1,2\dots p^{\lceil t \rceil}\}$, otherwise assume that $l(m_0)=l(m_1)=1$. According the \ref{karDef} definition $v_p(m_0-\alpha_i)\ge t$ and $v_p(m_1-\alpha_i)\ge t$. The triangle inequality implies that $v_p(m_0-m_1)\ge t$ which means that $m_0$ and $m_1$ in the same residue class modulo $p^{\lceil t \rceil}$. Therefore $m_0$ and $m_1$ must be equal if $m_0,m_1\in\{1,2\dots p^{\lceil t \rceil}\}$.
\\There exist $m_0\in \ZZ_p$ which maximizes $\min\{v_p(m-\alpha_i),v_p(m-\beta_j)\}$ otherwise $v_p(m-\alpha_i)$ and $v_p(m-\beta_j)$ would be arbitrary large for sufficient $m$ and due to the triangle inequality $v_p(\alpha_i-\beta_j)$ would be arbitrary large as well. However, the last one implies $\alpha_i=\beta_j$ which contradicts the assumption that the resultant is not zero.
\\Obviously,
$$\int_0^\infty \chi_t^{(x-\alpha_i)}(m_0)\chi_t^{(x-\beta_j)}(m_0)dt\le
\int_0^\infty \sum_{m=1}^{p^{\lceil t\rceil}}\kar{x-\alpha_i}\kar{x-\beta_j} dt$$
and since there is at most one non-zero term in each sum for a given $t$, we have
$$\int_0^T \chi_t^{(x-\alpha_i)}(m_0)\chi_t^{(x-\beta_j)}(m_0)dt\ge
\int_0^T \sum_{m=1}^{p^{\lceil t\rceil}}\kar{x-\alpha_i}\kar{x-\beta_j} dt$$
where $T=\min\{v_p(m_0-\alpha_i),v_p(m_0-\beta_j)\}$. On the other hand, both integrands are $0$ if $t>T$ otherwise $m_0$ would not satisfy the maximal condition. So, we got that
\begin{align*}
    \int_0^\infty \sum_{m=1}^{p^{\lceil t\rceil}}\kar{x-\alpha_i}\kar{x-\beta_j} dt&=\int_0^\infty \chi_t^{(x-\alpha_i)}(m_0)\chi_t^{(x-\beta_j)}(m_0)dt
    \\&=\max_{m\in\ZZ_p}\min\{v_p(m-\alpha_i),v_p(m-\beta_j)\}
\end{align*}
Finally, we need to put the pieces together.
\begin{align*}
    v_p(r)&\ge \sum_{i,j}\max_{m\in\ZZ_p}\min\{v_p(m-\alpha_i),v_p(m-\beta_j)\}\\
    &=\sum_{i,j} \int_0^\infty \sum_{m=1}^{p^{\lceil t\rceil}}\kar{x-\alpha_i}\kar{x-\beta_j} dt\\
    &=\int_0^\infty \sum_{m=1}^{p^{\lceil t\rceil}}\left(\sum_i\kar{x-\alpha_i}\right)\left(\sum_j\kar{x-\beta_j} \right)dt\\
    &=\int_0^\infty \sum_{m=1}^{p^{\lceil t\rceil}}\kar{f}\kar{g} dt
\end{align*}
\end{proof}
\begin{Def}
The integer transformed version of $\kar{h}$ is denoted by $\kkar{h}$ and defined by the equation:
$$\widehat{\chi}_r^{(h)}(m)=\int_{\lceil r\rceil-1}^{\lceil r\rceil}\kar{h} dt$$
\end{Def}
\begin{lemma}
The function, $\kkar{h}$ satisfies the following properties:
\begin{enumerate}
    \item It takes only integer values for every $m\in\ZZ_p$ and $t\in\RR_{\ge0}$.
    \item It is monotonically decreasing in $t$.
    \item $$v_p(h(m))=\int_0^\infty \kkar{h}dt=\sum_{t=1}^\infty \kkar{h}$$
\end{enumerate}
\end{lemma}
\begin{proof}
Let $h(x)=\prod_{i=1}^k (x-\delta_i)$ be a monic irreducible polynomial over $\QQ_p$ with $\delta_i$ roots in its splitting field $K$. Since $h$ is irreducible, we know that
$$v_p(m-\delta_i)=\frac{v_p(h(m))}{k}$$
Therefore $\kar{h}=k$ for every $t\in [0,\frac{v_p(h(m))}{k}]$ and $0$ otherwise. Since $v_p(h(m))$ is an integer, the length of the interval $[\lceil r \rceil -1, \lceil r\rceil]\cap [0, \frac{v_p(h(m))}{k}]$ can be written as $\frac{l}{k}$, where $l\in \ZZ$. So, we have that
\begin{align*}
    \widehat{\chi}_r^{(h)}(m)=\int_{\lceil r\rceil-1}^{\lceil r\rceil}\kar{h} dt=k\frac{l}{k}\in \ZZ
\end{align*}
On the other hand, $\kkar{h_1}+\kkar{h_2}=\kkar{h_1h_2}$ which implies that $\kkar{h}$ takes only integer values for every monic polynomial.
\\Since $\kar{h}$ is monotonically decreasing and we got $\kkar{h}$ by averaging out $\kar{h}$ on every $[\lceil r\rceil-1,\lceil r\rceil]$ interval, $\kkar{h}$ is also monotonically decreasing.
\\Finally, the last equation follows from the facts that $v_p(h(m))=\int_0^\infty \kar{h}dt=\int_0^\infty \kkar{h}dt$ and $\widehat{\chi}_r^{(h)}(m)=\int_{\lceil r\rceil-1}^{\lceil r\rceil}\kkar{h}$.

\end{proof}
\begin{lemma}\label{integralszummabecs}

\begin{align*}
    \int_0^\infty \sum_{m=1}^{p^{\lceil t\rceil}}\kar{f}\kar{g} dt&\ge \int_0^\infty \sum_{m=1}^{p^{\lceil t\rceil}}\kkar{f}\kkar{g} dt\\&=\sum_{t=1}^\infty \sum_{m=1}^{p^{t}}\kkar{f}\kkar{g} 
\end{align*}
\end{lemma}

\begin{proof}
Recall that $F_r:={\widehat{\chi}_r^{(f)}(m)}=\left(\int_{r-1}^{r}\kar{f}\right)$ and ${\widehat{\chi}_r^{(g)}(m)}=\left(\int_{r-1}^{r}\kar{g}\right)$. Since $\kar{f}$ and $\kar{g}$ are monotonically decreasing non-negative functions, the $\kar{f}-F_r$ is non-negative on the interval $[r-1, r_0]$ and negative on $(r_0, r]$. Furthermore, we have
$$U:=\int_{r-1}^{r_0} (\kar{f}-F_r)dt = - \int_{r_0}^r (\kar{f}-F_r)dt$$
where $U\ge 0$.
The monotonically decreasing property of $\kar{g}$ implies that
\begin{align*}
    &\int_{r-1}^{r_0} (\kar{f}-F_r)\kar{g} dt \ge U\chi_{r_0}^{(g)}(m) \ge -\int_{r_0}^r (\kar{f}-F_r)\kar{g}dt
\\\Rightarrow &\int_{r-1}^r (\kar{f}-F_r)\kar{g} dt \ge 0
\\\Rightarrow &\int_{r-1}^r \kar{f}\kar{g} dt \ge F_r\left(\int_{r}^{r+1}\kar{g}\right)={\widehat{\chi}_r^{(f)}(m)}{\widehat{\chi}_r^{(g)}(m)}
\end{align*}
The rest of the proof follows from the linearity aspect of $\kkar{.}$ and $\kar{.}$.

\end{proof}

\begin{lemma}
\label{divisonlemma}
\begin{align*}
    \widehat{\chi}_{t-1}^{(f)}(m)\ge \sum_{i=0}^{p-1}\widehat{\chi}_{t}^{(f)}(m+ip^{t-1})
\end{align*}
\end{lemma}

\begin{proof}
Since both sides of the inequality are additive, we only need to prove it for $f(x)=x-\delta$. Assume that $t\in(T-1,T]$ for some $T\in\ZZ$, then $\chi_t^{(x-\delta)}(m+ip^{T-1})$ is $1$ for at most one $i\in\{0,1\dots p-1\}$. Suppose the opposite, then $v_p(m+ip^{t-1}-\delta)\ge t$ and $v_p(m+jp^{t-1}-\delta)\ge t$, which implies that $v_p((i-j)p^{r-1})\ge t>r-1$, thus $v_p(i-j)> 0$, thus $i=j$. The $\kkar{x-\delta}$ is monotonically decreasing, and it only takes $0,1$ values, thus the sum is zero if the left side is zero, or at most one if the right side is one. 

\end{proof}

\section{Combinatorial Approach}
In the previous section, we managed to extract some beneficial equations. Although it has not been clarified yet, those equations let us handle the main problem in a combinatorial manner. In the subsequent chapter, we will show how to derive this form, but first, we want to lay the combinatorial foundation.

\begin{Def}
Let $\mathcal{F}=(V,E)$ be an infinite directed tree with $\delta$ root such as every vertex has $p$ out-degree.
\end{Def}
\begin{Def}
\label{abdef}
An $a:V\to \RR_{\ge1}\cup\{0\}$ weight function with $\omega_a$ weight is a function on the vertices of $\mathcal{F}$ satisfying the following properties.
\begin{enumerate}
    \item For every $T$ infinite path from $\delta$, $\sum_{v\in T}a(v)\ge \omega_a$.
    \item For every $v\in V$ vertex $a(v)\ge \sum_{u\in N_v} a(v)$, where $N_v$ denotes the set of end points of out-edges from $v$.
\end{enumerate}
Furthermore, if the range of $a$ is a subset of $\ZZ_{\ge 0}$ then we call $a$ an integral weight function and a real weight function otherwise.
\end{Def}
\begin{Def}
The scalar product of two weight is defined by the following equation:
$$\scal{a}{b}=\sum_{v\in V}a(v)b(v)$$
\end{Def}
Our goal is to minimize $\scal{a}{b}$ in terms of the weights of each weight functions, $\omega_a$ and $\omega_b$.
\begin{Def}
\label{gammadef}
We call a $(\gamma_i(\omega))_{i=0}^\infty\subset \RR_{\ge1}\cup\{0\}$ series the resolution of $\omega$ if
\begin{align}
    \gamma_i(\omega)&\ge p\gamma_{i+1}(\omega)\label{gammacsokken}\\
     \omega&=\sum_{i=0}^\infty \gamma_i(\omega)\label{gammaosszeg}
\end{align}
Then the lexicographically minimal is called the minimal resolution of $\omega$ and it is denote by $\rg_i$.
\\Furhetmore, a $(\gamma_i(\omega))_{i=0}^\infty\subset \ZZ_{\ge0}$ series is an integral resolution of $\omega$ and similarly, the lexicographically minimal is called the minimal integral resolution and is denoted by $\ig_i$.
\end{Def}
Note that $\omega\ge \gamma_i(\omega)$ and there exists and index $N$ such that $\gamma_i(\omega)=0$ for every $i\ge N$. Therefor the set of resolutions is compact and it has a lexicographically minimal element. Similarly, there exists a minimal integral resolution for every $\omega$.
\begin{lemma}
\label{realgammadef}
Let $k=\lfloor\log_p((p-1)\omega+1)\rfloor-1$, then the minimal resolution is
$$\rg_i(\omega)=p^{-i}\frac{p-1}{p-p^{-k}}\omega$$
where $i\le i\le k$ and $\rg_i=0$ otherwise.
\end{lemma}
\begin{proof}
First of all, we will show that there exist $x\in \RR$ and $k'$ such that $\rg_i=p^{-i}x$ for every $i\le k'$ and $\rg_i=0$ otherwise. Suppose the opposite, then there exists an $i$ index for which $\rg_i(\omega)> p\rg_{i+1}(\omega)\neq 0$. Let $\varepsilon>0$ be a sufficiently little real number such that $\rg_i(\omega)-\varepsilon\ge p(\rg_{i+1}(\omega)+\varepsilon)$. By replacing $\rg_i$ with $\rg_i-\varepsilon$ and $\rg_{i+1}$ with $\rg_{i+1}+\varepsilon$ we obtain a new real resolution of $\omega$ which differs form the initial one and is also lexicographically less. This contradicts our assumption, hence
\begin{align*}
    \rg_i(\omega)=\begin{cases}xp^{-i}&\text{if}~i\le k'\\0&\text{if}~i>k'\end{cases}
\end{align*}
According \ref{gammadef}, $\omega=x\sum_{i=0}^{k'} p^{-i}$ which allows us to determine $x$. On the other hand, $\rg_i\ge 1$ for every $i\le k'$, because the range of $\rg_i$ is $\RR_{\ge 1}\cap {0}$, thus we have $x\ge p^{k'}$ which implies that
\begin{align*}
    \omega=x\sum_{i=0}^{k'}p^{-i}=x\frac{p-p^{-k'}}{p-1}&\ge \frac{p^{k'+1}-1}{p-1}
    \\(p-1)\omega+1&\ge p^{k'+1}
    \\ \log_p((p-1)\omega+1)-1&\ge k'
    \\\lfloor\log_p((p-1)\omega+1)\rfloor-1&\ge k'
\end{align*}
In the last inequality we used the fact that $k'$ is an integer. We will show that $\rg_{k'}<p+\frac{p-1}{p^{k'+1}-1}$. Again, suppose the opposite and let $c=\frac{p-1}{p^{k'+1}-1}$. We define the following real resolution of $\omega$:
\begin{align*}
    \rg_i'=\begin{cases}
    \rg_i(\omega)-cp^{k'-i}&\text{for every}~ 0\le i\le k'
    \\c\sum_{j=0}^{k'}p^{k'-i}=1&\text{if}~i=k'+1
    \\0&\text{if}~ i>k'+1
    \end{cases}
\end{align*}
It is easy to check that $\rg_i'$ is indeed a resolution of $\omega$ and it is also lexicographically smaller than $\rg_i$ which contradicts our assumption. Hence, we have 
\begin{align*}
    \omega=x\frac{p-p^{-k'}}{p-1}&< \left(p+\frac{p-1}{p^{k'+1}-1}\right)\frac{p^{k'+1}-1}{p-1}=\frac{p^{k'+2}-1}{p-1}
    \\(p-1)\omega+1&< p^{k'+2}
    \\ \log_p\left((p-1)\omega+1\right)-2&< k'
\end{align*}
From the last two inequality it follows that $k'=\lfloor\log_p((p-1)\omega+1)\rfloor-1=k$ and $k=k'$. Finally, we can calculate each term of the series:
\begin{align*}
    \omega=x\frac{p-p^{-k}}{p-1}\implies x=\omega\frac{p-1}{p-p^{-k}}\implies \rg_i=p^{-i}\frac{p-1}{p-p^{-k}}\omega
\end{align*}
where $i\le k$. Note that the definition of $k'$ implies $\rg_i=0$ if $i>k$.
\end{proof}
\begin{lemma}
\label{gammalemma}
Let $\gamma$ be either a real or integral minimal resolution, then $\gamma$ satisfies the following properties.
\begin{enumerate}
    \item $\gamma_1(\omega)=\gamma_0(\omega-\gamma_0(\omega))$.
    \item $\gamma_i(\omega)$ monotonically increasing in $\omega$.
\end{enumerate}
\end{lemma}
\begin{proof}
Note that the series $(\gamma_{i+1}(\omega))_{i=0}^\infty$ is a minimal resolution with $\omega-\gamma_0(\omega)$ weight, therefore $\gamma_1(\omega)=\gamma_0(\omega-\gamma_0(\omega))$.

For the second property suppose the opposite. More precisely, let $\omega_1>\omega_2$ be two weight such that $\gamma_0(\omega_1)<\gamma_0(\omega_2)$. By reducing $\gamma(\omega_1)$ to a $\gamma'$ resolution with $\omega_2$ weight, we obtain a lexicographically smaller resolution than $(\gamma_i(\omega_2))_{i=0}^\infty$, which contradicts the fact that $\omega_i$ is a minimal resolution. Let $k$ be the smallest non-negative integer such that $\sum_{i=0}^k \gamma_i(\omega_1) \ge \omega_2$, then we can construct $\gamma'$ as
\begin{align*}
    \gamma_i' = \begin{cases}
    \gamma_i(\omega_1)~&\text{if}~i<k
    \\\omega_2-\sum_{j=0}^{i-1}\gamma_i(\omega_1)~&\text{if}~i=k
    \\0~&\text{if}~ i>k
    \end{cases}
\end{align*}
Furthermore, assume that $\omega_1-\omega_2<\gamma_0(\omega_1)-\gamma_0(\omega_2)$, then we can define the following resolution
$$\gamma_i'=\begin{cases}\gamma_0(\omega_2)+\omega_1-\omega_2&\text{if}~i=0\\\gamma_i(\omega_2)&\text{if}~i>0\end{cases}$$
Due to our assumption, $\gamma_i'$ is lexicographically smaller than $\gamma_i(\omega_1)$, which is a contradiction. So, we showed that both $\gamma_0(\omega)$ and $\omega-\gamma_0(\omega)$ are monotonically increasing in $\omega$. 
Let $\omega^{(k)}$ denote the weight of $(\gamma_i(\omega))_{i=k+1}^\infty$. Note that $\omega^{(0)} = \omega-\gamma_0(\omega)$. We will prove by induction that $\omega^{(k)}$ and $\gamma_k(\omega)$ are monotonically increasing in $\omega$. First, we have already established the base case, when $k=0$. Then, assume that the statement holds for every $k<K$. The $(\gamma_i(\omega))_{i=K}^\infty$ is also a minimal resolution, thus $\gamma_0(\omega^{(K-1)}) =\gamma_K(\omega)$ and $\omega^{(K)} = \omega^{(K-1)}-\gamma_0(\omega^{(K-1)})$. We know that $\gamma_0(\omega')$ and $\omega'-\gamma_0(\omega')$ are monotonically increasing in $\omega'$ and $\omega^{(K-1)}$ is also monotonically increasing in $\omega$, consequently $\gamma_0(\omega^{(K-1)})=\gamma_K(\omega)$ and $\omega^{(K-1)}-\gamma_0(\omega^{(K-1)}) = \omega^{(K)}$ are monotonically increasing. Hence we proved the induction step.
\end{proof}
\begin{lemma}
\label{forind}
Let $\gamma$ be an integral or a real minimal resolution and $\omega<\sum_{i=0}^k p^i$, then we have $\omega-\gamma_0(\omega)<\sum_{i=0}^{k-1} p^i$.
\end{lemma}
\begin{proof}
Let $k$ be the minimal integer that satisfies the condition of the lemma, then $k$ is the same defined in \ref{realgammadef}. First, we prove the lemma for real minimal $\gamma$. Due to \ref{realgammadef}, we have $\rg_0(\omega)=\frac{p-1}{p-p^{-k}}\omega$. By using $\omega<\sum_{i=0}^k p^i$, we get
\begin{align*}
    \omega-\rg_0(\omega)<\left(1-\frac{p-1}{p-p^{-k}}\right)\sum_{i=0}^{k} p^i=\frac{1-p^{-k}}{p-p^{-k}}\frac{p^{k+1}-1}{p-1}=\sum_{i=0}^{k-1}p^i
\end{align*}
We proved the real case. Now, we can move on to the integral case. Since every integral resolution is also a real resolution $\rg$ has to be lexicographically smaller than $\ig$. The latter one implies that $\ig_0(\omega)\ge \rg_0(\omega)\Rightarrow \omega-\ig_0(\omega)\le \omega-\rg_0(\omega)$.
\end{proof}

\begin{lemma}
\label{rootest}
Let $a$ be a real (or integral) weight function with $\omega_a$ weight and $\gamma$ real (or integral) minimal resolution of $\omega_0$. Then $a(\delta)\ge \gamma_0(\omega_a)$.

\end{lemma}

\begin{proof}
We start with the real case. For every $v\in V$ there exists an $u\in N_v$ such that $a(u)\le \frac{a(v)}{p}$, otherwise it follows from the definition that $a(v)\ge \sum_{u\in N_v} a(u)> |N_v|\frac{a(v)}{p}=p\frac{a(v)}{p}$, which is clearly a contradiction. Therefore, we have an infinity path from the root $T=(\delta=u_0, u_1, \dots)$ with the property that for every $u_k\in T$
$$a(u_k)\le \frac{a(u_{k-1})}{p}\le \dots \le \frac{u(\delta)}{p^k}$$
Let $k'$ denote the minimum index such that $a(u_{k'+1})=0$. From the definition of real-valued weight functions, we have
$$\omega_a \le \sum_{u\in T} a(u)\le \sum_{i=0}^{k'} \frac{a(\delta)}{p^i}$$
Let us define the $\gamma_i'=\frac{a(\delta)}{p^i}$ series, which is clearly a resolution. If $\omega'$ is the weight of $\gamma'$, then the last inequality means that $\omega_a\le \omega'$. Since $\rg(\omega')$ lexicographically smaller than $\gamma'$ (Because $\ig(\omega')$ was defined to be minimal.), we also have $a(\delta)=\gamma_0'\ge\ig_0(\omega')$. Finally, \ref{gammalemma} implies $a(\delta)\ge \rg_0(\omega')\ge \rg_0(\omega_a)$ because $\omega'\ge \omega_a$.

Now, we move on to the integral case. Due to the fact that every element is an integer, we have a slightly stronger inequality, namely $a(u)\le \left\lfloor\frac{a(\delta)}{p^k}\right\rfloor$. Similarly to the real case, this implies
$$\omega_a\le \sum_{v\in T}a(v)\le \sum_{i=0}^\infty \left\lfloor\frac{a(\delta)}{p^k}\right\rfloor$$
Once again lets take $\gamma_i'=\left\lfloor\frac{a(\delta)}{p^k}\right\rfloor$. It is easy to see that $\gamma_i'$ is an integral resolution and therefore $\omega'\ge \omega_a$ (where $\omega'$ is the weight of $\gamma'$). Since $\ig_0(\omega')$ is lexicographically minimal and monotonically increasing in $\omega'$, we get that $a(\delta)=\gamma_0\ge \ig_0(\omega')\ge \ig_0(\omega_a)$.
\end{proof}
\begin{Th}
\label{mainoptth}
Let $a,b$ be real (or integral) weight functions and $\gamma$ the real (or inegral) minimal resolution, then we have the following inequality
$$\scal{a}{b}\ge \sum_{i=0}^{\infty}p^i\gamma_i(\omega_a)\gamma_i(\omega_b)$$
\end{Th}
\begin{proof}
Assume that $\omega_a<\sum_{i=0}^k p^i$, $\omega_a\le\omega_b$, where $k\in\ZZ_{>0}$. We apply induction on $k$. First, we show the base case, when $k=1$. It follows from \ref{forind} that $\gamma_0(\omega_a)>\omega_a-1$. We also know that $\gamma_i(\omega_a)=0$ if $i\ge 1$, otherwise $\gamma_j(\omega_a)\ge$ for some $j\ge 1$, and $\omega_a =\sum_{i=0}^\infty \gamma_i(\omega_a) \ge \gamma_0(\omega) +\gamma_j(\omega_a)>\omega_a$, which is clearly contradiction.
Since every term is non-negative, we have $\scal{a}{b}\ge a(\delta)b(\delta)$. By using \ref{rootest} and the fact that $\gamma_i(\omega_a)=0$ if $i\ge 1$, we obtain the following inequality:
$$\scal{a}{b}\ge a(\delta)b(\delta)\ge \gamma_0(\omega_a)\gamma_0(\omega_b)=\sum_{i=0}^{\infty}p^i\gamma_i(\omega_a)\gamma_i(\omega_b) \ .$$
Hence, we proved the theorem for $k=1$.

Let $k>1$ and assume that the statement is true for every integer smaller than $k$. The \ref{forind} lemma states that $\omega_a-\gamma_0(\omega_a)<\sum_{i=0}^{k-1}p^i$. Let $v$ be an arbitrary vertex of the tree ($v\in V$) and let $a^v$ denote the the constraint of the weight function to a sub-tree from the root $v$. (Every ancestor of $v$ is included in the $a^v$ graph.) Moreover, let $\omega_a^v$ denote the weight of $a^v$. This new notion allows us to rewrite the scalar product
$$\scal{a}{b}= a(\delta)b(\delta)+\sum_{v\in N_\delta}\scal{a^v}{b^v}$$
According the induction hypothesis, we already proved the statement for $\scal{a^v}{b^v}$ because $\omega^v_a=\omega_a-\gamma_0(\omega_a)<\sum_{i=0}^{k-1}p^i$. Hence, we obtain the following inequality

\begin{align*}
\scal{a}{b}&\ge a(\delta)b(\delta)+\sum_{v\in N_\delta}\sum_{i=0}^{\infty}p^i\gamma_i(\omega^v_a)\gamma_i(\omega^v_b)\\&=a(\delta)b(\delta)+p\sum_{i=0}^{\infty}p^i\gamma_i(\omega^v_a)\gamma_i(\omega^v_b)
\end{align*}

Recall that $\omega_a^v=\omega_a-a(\delta)$, therefore the right side of the last inequality only depends on $a(\delta)$ and $b(\delta)$. We claim that the right side is minimal if $b(\delta)=\rg_0(\omega_b)$. For the sake of contradiction, assume that $b(\delta)>\rg_0(\omega_b)$. (Recall that \ref{rootest} states that $b(\delta)\ge\rg_0(\omega_b)$.) Let $\varepsilon=b(\delta)-\gamma_0(\omega_b)$. We want to prove that the right side will not increase if we decrease $b(\delta)$ by $\varepsilon$. 

\begin{align*}
    a(\delta)b(\delta)+p\sum_{i=0}^{\infty}p^i\gamma_i(\omega^v_a)\gamma_i(\omega^v_b) \ge a(\delta)(b(\delta)-\varepsilon)+p\sum_{i=0}^{\infty}p^i\gamma_i(\omega^v_a)\gamma_i(\omega^v_b+\varepsilon)
\end{align*}
Or equivalently, 
\begin{align*}
    a(\delta)\varepsilon \ge p\sum_{i=0}^\infty p^i\gamma_i(\omega_a^v)\left(\gamma_i(\omega_b^v+\varepsilon)-\gamma_i(\omega_b^v)\right)
\end{align*}

Let $\Delta_i=\gamma_i(\omega_b^v+\varepsilon)-\gamma_i(\omega_b^v)$ denote the difference between each term. (Recall that $\Delta_i\ge 0$ due to \ref{gammalemma} lemma.) Also the \ref{gammadef} definition implies that $\sum \Delta_i=(\omega_b^v+\varepsilon)-\omega_b^v=\varepsilon$. Therefore it is enough to prove that
\begin{align*}
    a(\delta)\Delta_i\ge p^{i+1}\gamma_i(\omega_a^v)\Delta_i
    \Leftarrow a(\delta)\ge p^{i+1}\gamma_i(\omega_a^v)
\end{align*}
Because we would obtain the desired inequality if we sum them up. Fortunately, these inequalities can be proved with the following chain of inequalities: 
\begin{align*}
    a(\delta)&\overset{(1)}{\ge} \gamma_0(\omega_a)\overset{(2)}{\ge} p\gamma_1(\omega_a)\overset{(3)}{=}p\gamma_0(\omega_a-\gamma_0(\omega_0))
    \\&\overset{(4)}{\ge} p\gamma_0(\omega_a-a(\delta))=p\gamma_0(\omega_a^v)\overset{(2)}{\ge} p^{i+1}\gamma_i(\omega_a^v)
\end{align*}
Where (1) is \ref{rootest}, (2) is \ref{gammadef}.1, (3) is \ref{gammalemma}.1 and (4) follows from \ref{gammalemma}.2 because $a(\delta)\ge \gamma_0(\omega_a)$ due to $\ref{rootest}$. We completed the proof.
\end{proof}
\begin{Cor}
\label{valosegeszsulyfv}
Every real weight function is also integral and therefore we have the following inequality for the $a,b$ integral weight functions
\begin{align*}
    \scal{a}{b}\ge \sum_{i=0}^{\infty}p^i\ig_i(\omega_a)\ig_i(\omega_b)\ge \sum_{i=0}^{\infty}p^i\rg_i(\omega_a)\rg_i(\omega_b)
\end{align*}
\end{Cor}
\begin{Cor}
\label{computedrealgamma}
For arbitrary $a$ and $b$ real weight function let $k_c=\lfloor\log_p((p-1)\omega_c+1)\rfloor-1$ for every $c\in\{a,b\}$. Then we have
\begin{align*}
    \scal{a}{b}\ge\sum_{i=0}^{\infty}p^i\rg_i(\omega_a)\rg_i(\omega_b)= \frac{p-1}{p-p^{-\max\{k_a,k_b\}}}\omega_a\omega_b
\end{align*}
\end{Cor}
\begin{proof}
Due to \ref{realgammadef}, we can substitute the values of $\rg_i$ into \ref{mainoptth} theorem, which gives us
\begin{align*}
    \scal{a}{b}\ge \frac{p-1}{p-p^{-k_a}}\frac{p-1}{p-p^{-k_b}}\omega_a\omega_b\sum_{i=0}^{\min\{k_a,k_b\}}p^{-i}
\end{align*}
where the sum goes only to $\min\{k_a,k_b\}$ because during the proof of \ref{realgammadef} we showed that $\rg_i(\omega_a)=0$ for every $k$ index greater than $k_a$. Lastly,
\begin{align*}
    \sum_{i=0}^{\min\{k_a,k_b\}}p^{-i} =\frac{p-p^{-\min\{k_a,k_b\}}}{p-1}
\end{align*}
Substituting it back yields the inequality.
\end{proof}
\begin{lemma}
\label{sulyfvbecslese}
Let $\gamma$ be a minimal resolution (real or integral) and introduce $a(v)=\gamma_{d(\delta,v)}(\omega_a)$ weight function with $\omega_a$ weight, where $d(\delta, v)$ denotes the distance of the $v$ vertex from the $\delta$ root. Similarly, we define the $b$ weight function with $\omega_b$ weight. Then the lower bound in \ref{mainoptth} theorem is sharp.
\end{lemma}
\begin{proof}
Since there are $p^i$ vertices $i$ distance from the $\delta$ root. With little calculation we get that
\begin{align*}
    \scal{a}{b}=\sum_{v\in V}a(v)b(v)=\sum_{i=0}^\infty p^i\gamma_i(\omega_a)\gamma_i(\omega_b)
\end{align*}
which is exactly the lower bound in \ref{mainoptth} theorem. In conclusion the provided inequality in \ref{mainoptth} is sharp.
\end{proof}
\section{Proof of Main Theorem}
In the last two section we prepared every tool we will need in the proof of the main theorem. In this chapter we will focus on that proof and the corollaries of the main theorem.
\begin{Th}
\label{maintheorem}
Let $f,g$ be monic polynomials with integer coefficients and $r$ resultant such that $v_p(f(n))\ge s_1$ and $v_p(g(n))\ge s_2$ for every integer $n$. Furthermore, let $\gamma$ be the (real or integral) minimal resolution defined in \ref{gammadef}. Then we have the following lower bound for the $p$-adic valuation of the $r$ resultant.
\begin{align*}
    v_p(r)\ge p\sum_{i=0}^{\infty}p^i\gamma_i(s_1)\gamma_i(s_2)
\end{align*}
\end{Th}
\begin{proof}
When $r=0$, $v_p(0)=\infty$, thus we can assume that $r\neq 0$.
Recall that the \ref{vprchiint} and \ref{integralszummabecs} lemmas imply the following inequality
\begin{align}
\label{mainsth}
    v_p(r)\ge \sum_{t=1}^\infty \sum_{m=1}^{p^t}\kkar{f}\kkar{g}
\end{align}
Take $p$ instances of the graph defined in the 3rd section and let $\mathcal{F}_0,\mathcal{F}_1,\dots,\mathcal{F}_{p-1}$ denote these instances. We identify the vertices of $\mathcal{F}_k$ by the $(l,m)$ pairs ($l\in \ZZ_{\ge1}$) with the following rule. 
Let $\varphi_i:\ZZ/p^{i+1}\ZZ\rightarrow\ZZ/p^i\ZZ$ be the natural projection and note that $k\xrightarrow{\varphi_1^{-1}}D_1\xrightarrow{\varphi_2^{-1}}D_2\rightarrow\dots$ determines a similar graph structure as $\mathcal{F}_k$, where $D_i=\{x\in \ZZ/p^{i+1}\ZZ:x\equiv k~(\text{mod}~p)\}$. 
More precisely, if $v\in \mathcal{F}_k$ then we identified $v$ with the pair $(l, m)$, where $m\in D_{d(\delta,v)}$, $l=d(\delta, v)$ and the end points of out-edges are $N_v=\{(l+1, m') \mid m'\in D_{l+1},~ \varphi_{l+1}(m')=m\} $. (Recall that $d(\delta,v)$ is the distance between $\delta$ and $v$. Furthermore, it is easy to see that $|N_v|=p$.)

We claim that $a(t,m)=\kkar{f}$ is an (integral) weight function. First of all, the weight of $a$ is $s_1$ because
$$\sum_{t=1}^\infty \kkar{f}=v_p(f(m))\ge s_1$$
for arbitrary $m\in\ZZ$. We also need to mention that due to $\ref{divisonlemma}$ lemma, the \ref{abdef}.2 criteria is also satisfied.  Similarly, $b(t,m)=\kkar{g}$ is weight function with $s_2$ weight.

Notice that the scalar product of $a$ and $b$ on $\mathcal{F}_k$ is
\begin{align*}
    \scal{a}{b}=\sum_{v\in \mathcal{F}_k} a(v)b(v)=\sum_{t=1}^\infty \sum_{\substack{m\equiv k~(p)\\ m\in\ZZ/p^{t+1}}}\kkar{f}\kkar{g}
\end{align*}
If we add these up for $0\le k<p$, we get the right side of \ref{mainsth}. Therefore, within in the meaning of \ref{mainoptth} theorem we get that
\begin{align*}
    v_p(r)\ge p\sum_{i=0}^{\infty}p^i\gamma_i(\omega_a)\gamma_i(\omega_b)\ .
\end{align*}
Hence, we finished the proof.
\end{proof}
Initially, the authors of the \cite{frenkel2018estimating} paper considered a problem setting, where $s_1=s_2$ and $S=\max_{n\in\ZZ}\{\min{v_p(f(n)), v_p(g(n))}\}$ was also known. This encouraged us to investigate this particular case as well.
\begin{Th}
Let $f,g$ be monic polynomials with integer coefficients and $r$ resultant such that $v_p(f(n))\ge s_1$, $v_p(g(n))\ge s_2$ for every integer $n$ and 
$$S=\max_{n\in\ZZ}\min \{v_p(f(n)),v_p(g(n))\}<\infty$$
Let $\gamma$ be the (real or integral) minimal resolution defined in \ref{gammadef}. Then we have
$$v_p(r)\ge S-\max\{s_1, s_2\}+p\sum_{i=0}p^i\gamma_i(s_1)\gamma_i(s_2)$$
\end{Th}
\begin{proof}
First, we need to understand what the statement on $S$ is really about. Since $S$ is not infinity and $\min \{v_p(f(n)),v_p(g(n))\}$ takes only non-negative integer values, there exists an $m_0\in\ZZ_p$ such that $S=\min \{v_p(f(m_0)),v_p(g(m_0))\}$. Similarly to the proof of theorem \ref{maintheorem}, let $\mathcal{F}_0,\dots \mathcal{F}_{p-1}$ denote the trees and $a,b$ the corresponding weight functions on them derived from $\kar{f}$ and $\kar{g}$.
Let $k\equiv m_0~(\text{mod}~p)$ and $\mathcal{T}$ denote the path in $\mathcal{F}_k$ determined by $m_0$. We can look up the values of the weight functions, $a,b$, along the $\mathcal{T}$ path. Let these number be $\alpha_0, \alpha_1, \dots $ and $\beta_0, \beta_1, ...$. (In other words, $\alpha_i=a(t_i)$, where $\mathcal{T}=(t_0,t_1,\dots)$.)
We already know that $v_p(f(m_0)), v_p(g(m_0))\ge S$, which implies $\sum_i \alpha_i, \sum_i \beta_i \ge S$. We can form a new $\alpha^*_0, \dots \alpha^*_u$ finite series from $(\alpha_i)_{i\in\NN}$ which satisfies $\alpha^*_i=\alpha_j$ if $j<u$ and $\sum_{i=0}^k \alpha^*_i=s_1$. (Where we used the fact that $s_1\le S$.) We can do the same thing for the series $(\beta_i)_{i\in\NN}$ and obtain $\beta^*_0, \dots \beta^*_v$.

By changing the values of $a$ along the $\mathcal{T}$ path to $\alpha^*_i$ we get a new weight function $a^*$. We claim that it also has weight $s_1$. The first changed value is $\alpha_u$ and every path containing that node already has $\sum_{i=0}^u \alpha^*_i=s_1$ weight while the other path remained untouched. So, it has at least $s_1$ weight. We should also check the 2nd constraint in Definition \ref{abdef}. Unfortunately, it is not necessary true; however, if we set $a^*(v)$ to $0$ for every $v$ ancestor of $t_u$ (where $\alpha_u$ is located), the weight of $a^*$ is still $s_1$ and the 2nd constraint also holds. Again, we can do the same thing with $b$ and obtain $b^*$. 

Let $\mu,\nu$ denote the indexes of the first zeros in $\alpha$ and $\beta$. (Possibly infinity.) If $\mu\le \nu$, then we have the following estimate.
\begin{align*}
    \scal{a}{b}&= \scal{a^*}{b} +\sum_{i=0}^\mu (\alpha_i-\alpha^*_i)\beta_i
    \\&\ge \scal{a^*}{b^*}+S-s_1
\end{align*}
Where we used that  $b^*\le b \Rightarrow \scal{a^*}{b^*}\le \scal{a^*}{b}$ and $\mu\le \nu$ implying that $\beta_i\ge 1$ for every $i\le \mu$, which means that 
$\sum_{i=0}^\mu (\alpha_i-\alpha^*_i)\beta_i\ge \sum_{i=0}^\mu (\alpha_i-\alpha^*_i) \ge S-s_1$. 
We can do the same thing if $\mu\ge \nu$. Combining the two together, we get that
$$\scal{a}{b}\ge \scal{a^*}{b^*} +\min\{S-s_1,S-s_2\}$$
Since $a^*,b^*$ are weight functions with $s_1, s_2$ weights, we can apply Theorem \ref{mainoptth} and get
$$\scal{a}{b}\ge  S-\max\{s_1, s_2\}+\sum_{i=0}p^i\gamma_i(s)$$

In the other trees we don't have any new constraint compared to Theorem \ref{maintheorem}, thus our final lower bound is 
$$v_p(r)\ge S-\max\{s_1, s_2\}+p\sum_{i=0}p^i\gamma_i(s)$$
\end{proof}
If we substitute our computations from \ref{computedrealgamma} when $\gamma$ is the real minimal resolution we get the following corollaries

\begin{Cor}
Let $f,g$ be monic polynomials with integer coefficients and $r$ resultant such that $v_p(f(n))\ge s_1$, $v_p(g(n))\ge s_2$ for every integer $n$ and 
$$S=\max_{n\in\ZZ}\min \{v_p(f(n)),v_p(g(n))\}<\infty$$
Then we have the following lower bound for $v_p(r)$.
\begin{align*}
    v_p(r)\ge S-\max\{s_1,s_2\}+ ps_1s_2\frac{p-1}{p-p^{-k}}
\end{align*}
Where $k=\lfloor \log_p((p-1)\max\{s_1,s_2\}+1)\rfloor -1$.
\end{Cor}
\begin{Cor}
Let $f,g$ be monic polynomials with integer coefficients and $r$ resultant such that $v_p(f(n))\ge s$, $v_p(g(n))\ge s$ for every integer $n$ and 
$$S=\max_{n\in\ZZ}\min \{v_p(f(n)),v_p(g(n))\}<\infty$$
We have
\begin{align*}
   v_p(r)\ge S-s+ ps^2\frac{p-1}{p-p^{1-\lfloor \log_p((p-1)s+1)\rfloor}}
\end{align*}
\end{Cor}

\section{Constructions}

In this section we discuss the sharpness of the provided lower bounds from the previous section. We show that the estimations from the $\gamma$ real minimal resolution give us the best lower bound infinitely many times.

\begin{lemma}
Let $\log_p((p-1)s_1+1)-1, k_2=\log_p((p-1)s_2+1)-1$ integers and $s_1\ge s_2 \Leftrightarrow k_1\ge k_2$. We define $f,g$ polynomials to be
\begin{align*}
    f(x)=\prod_{t=0}^{p-1}h(x+t)
    \\g(x)=\prod_{t=0}^{p^{k_2+1}-1}(x+t)
\end{align*}
where $h$ is a monic irreducible polynomials over $\QQ_p$ with degree $s_1$ such that $v_p(h(0))=s_1$ and $v_p(h(n))=0$ for every other $n$ that is not divisible by $p$. This two polynomial satisfy $v_p(f(n))\ge s_1$ and $v_p(g(n))\ge s_2$ for every $n$ integer and their resultant is
\begin{align*}
    v_p(r)=ps_1s_2\frac{p-1}{p-p^{-k_1}}
\end{align*}
(Which is the provided lower bound in Theorem \ref{maintheorem}.)
\end{lemma}
\begin{proof}
First of all, we prove the existence of the polynomial $h$. It is well known that there exists irreducible polynomial over $\mathbb{F}_p$ with degree $s_1$. Let $h_0$ be such polynomial. We claim that $h=p^{s_1}h_0(\frac{x}{p})$ is good choice. It was irreducible over $\mathbb{F}_p$, so it is also irreducible over $\QQ_p$. On the other hand, $v_p(h(0))=s_1$ because $h_0(0)\equiv 0~(\text{mod}~p)$ would mean that $h_0$ is not irreducible. Meanwhile, $h(x)\equiv cx^{s_1} ~(\text{mod}~p)$, where $c$ is the non-zero leading coefficient of $h_0$. Therefore, $v_p(h(n))=0$ for every other $n$ that is not divisible by zero.

We can factor $h$ over its splitting field and write $h=\prod_{i=1}^{s_1}(x-\alpha_i)$. Due to the irreducibility of $h$, we have $v_p(n-\alpha_i)=\frac{1}{s_1}v_p(h(n))$ or equivalently \begin{align*}
    v_p(n-\alpha_i)=\begin{cases}1&\text{if}~p|n\\ 0&\text{otherwise}\end{cases}
\end{align*}
Hence, the resultant is
\begin{align*}
    v_p(r)=\sum_{t=0}^{p^{k_2+1}-1}\sum_{i=1}^{s_1}v_p(t-\alpha_i)=p^{k_2+1}s_1
\end{align*}
but $s_2=\frac{p^{k_2+1}-1}{p-1}$ so
\begin{align*}
    v_p(r)=ps_1s_2p^{k_2}\frac{p-1}{p^{k_2+1}-1}=ps_1s_2\frac{p-1}{p-p^{-k_2}}
\end{align*}
Last of all, we need to check that $v_p(f(n))\ge s_1$ and $v_p(g(n))\ge s_2$ holds for every $n\in \ZZ$. 
\\In the case of $g$, we already showed that for an arbitrary $n\in\ZZ$, we have
$$v_p(g(n))=\sum_{t=0}^{p-1}\sum_{i=1}^{s_1} v_p(n+t-\alpha_i)$$
There exists $t_0\in\{0,1\dots p-1\}$ such that $p| n+t$. As we had showed earlier $v_p(n+t-\alpha_i)=1$, which means $v_p(g(n))=s_1$.
\\In the case of $f$, for a given $n\in \ZZ$ there are precisely $p^{k_2}$ integers $t\in[0,p^{k_2+1}-1]$ such that $v_p(n+t)\ge 1$. In general, there are $p^{k_2-l}$ integers $t$ for which $v_p(n+t)\ge l+1$. This gives us that 
\begin{align*}
    v_p(f(n))&=\sum_{l=0}^{k_2}p^{k_2-l}=\frac{p^{k_2+1}-1}{p-1}
    \\&=\frac{p^{\log_p((p-1)s_2+1)}-1}{p-1}=s_2
\end{align*}
where we substituted the definition of $k_2=\log_p((p-1)s_2+1)-1$.
We proved the lemma.
\end{proof}

\bibliographystyle{unsrt}
\bibliography{main.bib}

\end{document}